# Some Two-Step Procedures for Variable Selection in High-Dimensional Linear Regression

Jian Zhang†,  Xinge Jessie Jeng†,   Han Liu‡

† Department of Statistics,    ‡Department of Statistics

Purdue University    Carnegie Mellon University

October 22, 2018

## ABSTRACT

We study the problem of high-dimensional variable selection via some two-step procedures. First we show that given some good initial estimator which is $\ell_\infty$-consistent but not necessarily variable selection consistent, we can apply the nonnegative Garrote, adaptive Lasso or hard-thresholding procedure to obtain a final estimator that is both estimation and variable selection consistent. Unlike the Lasso, our results do not require the irrepresentable condition which could fail easily even for moderate $p_n$ (Zhao and Yu, 2007) and it also allows $p_n$ to grow almost as fast as $\exp(n)$ (for hard-thresholding there is no restriction on $p_n$). We also study the conditions under which the Ridge regression can be used as an initial estimator. We show that under a relaxed identifiable condition, the Ridge estimator is $\ell_\infty$-consistent. Such a condition is usually satisfied when $p_n \leq n$ and does not require the partial orthogonality between relevant and irrelevant covariates which is needed for the univariate regression in (Huang et al., 2008). Our numerical studies show that when using the Lasso or Ridge as initial estimator, the two-step procedures have a higher sparsity recovery rate than the Lasso or adaptive Lasso with univariate regression used in (Huang et al., 2008).

## I. Introduction

Consider the linear regression model

$$Y = X\beta^* + \epsilon \qquad (1.1)$$

where $X \in \mathbb{R}^{n \times p}$ is the design matrix, $Y \in \mathbb{R}^{n \times 1}$ is the response vector, $\beta^* \in \mathbb{R}^{p \times 1}$ is the unknown parameter, and errors $\epsilon = [\epsilon_1, \ldots, \epsilon_n]^T$ are iid normal, i.e. $\epsilon \sim N(0, \sigma^2 I)$. We are interested in regression with diverging number of parameters, and will use $p_n$ to denote the number of variables which can grow as $n \to \infty$.

The key assumption for such high-dimensional estimation problems to be feasible is that the true parameter $\beta^*$ is sparse. Let $S$ be the subset of indices such that $S = \{j | \beta_j^* \neq 0\}$ and denote $s_n = |S|$, the cardinality of the set $S$. The sparsity assumption means that the number of relevant variables $s_n$ is much smaller than $p_n$, i.e. $s_n \ll p_n$. Under such a condition, efficient estimation and variable selection become possible. For example, the Lasso (Tibshirani, 1996) which minimizes least squares with the $\ell_1$ penalty

$$\beta^{Lasso} = \arg\min \frac{1}{2n} \|Y - X\beta\|^2 + \lambda_n \sum_{j=1}^{p_n} |\beta_j| \qquad (1.2)$$

has been proposed for such problems. Due to the $\ell_1$ penalty, the solution of Lasso is usually sparse with an appropriately chosen penalty parameter $\lambda_n$. Such a property has made Lasso a very desirable candidate for variable selection. Computationally, the estimation of Lasso is a convex optimization problem and can be solved efficiently. Furthermore, it has been shown that the full solution path of Lasso can be found at the same cost of solving the least squares estimation problem (Osborne et al., 2000; Efron et al., 2004). People have also studied various theoretical properties of Lasso (Fu and Knight, 2000; Greenshtein and Ritov, 2004; Meinshausen and Bühlmann, 2006; Zou, 2006; Zhao and Yu, 2007; Yuan and Lin, 2007; Bickel et al., 2007; Wainwright, 2006). One interesting property found by several authors (Meinshausen and Bühlmann, 2006; Zou, 2006; Zhao and Yu, 2007) is that Lasso is not variable selection consistent in general, and a condition on the design matrix (called the irrepresentable condition in (Zhao and Yu, 2007)) is needed to ensure its variable selection consistency. For high-dimensional inference with increasing $p_n$, several studies (Meinshausen and Bühlmann, 2006; Zhao and Yu, 2007; Wainwright, 2006) showed that under the irrepresentable condition, Lasso is also variable selection consistent if additional conditions on $p_n$, $s_n$, $n$ and $\lambda_n$ are satisfied. In particular, it has been shown that $p_n$ can be allowed to grow almost as fast as $\exp(n)$ when the error is normally distributed. Although such theoretical results are very encouraging for the Lasso in high-dimensional problems, it has been pointed out in (Zhao and Yu, 2007) that the key irrepresentable condition on the design matrix can easily fail even for moderate $p_n$.

On the other hand, it is shown in (Fan and Li, 2001; Zou, 2006) that even if the irrepresentable condition is satisfied and the Lasso is variable selection consistent, there does not exist a tuning parameter $\lambda_n$ which can lead to both efficient estimation and consistent variable selection. It is argued that the desired estimator should possess the oracle properties (Fan and Li, 2001), i.e. it should be variable selection consistent and the estimation of the nonzero parameters should be efficient. As a result, the SCAD method has been proposed and studied for both the fixed and increasing $p_n$ setting with $p_n^5/n \to 0$ (Fan and Li, 2001; Fan and Peng, 2004), and it has been shown to have the oracle properties. Huang et al. (2008) showed that the bridge estimator (Frank and Friedman, 1993) for linear model, which has a penalty term $\lambda_n \sum_{j=1}^{p_n} |\beta_j|^\gamma$ for $0 < \gamma < 1$, also has the oracle properties under certain conditions when $p_n < n$. However, since the penalty functions of both the SCAD and the bridge estimator are non-convex, it is more difficult to solve such optimization problems and in general there is no guarantee to find the global minimizer efficiently especially when the number of variables is large.

Recently several two-step procedures have been studied for variable selection. The adaptive Lasso approach, which was recently proposed by Zou (2006), uses a weighted $\ell_1$ penalty with weights determined by an initial estimator. In other words, the adaptive Lasso can be thought as a two-step procedure by applying the Lasso to some transformed design with the initial estimator. For fixed $p_n$, Zou (2006) showed that if the initial estimator satisfies certain conditions related to estimation consistency, the adaptive Lasso estimator has the oracle properties. Huang et al. (2006) further extended the results of the adaptive Lasso with increasing $p_n$. Yuan and Lin (2007) studied the nonnegative Garrote method (Breiman, 1995) for fixed $p_n$ and proved that when supplied with some good initial estimator which is $\ell_\infty$-consistent, the final nonnegative Garrote estimator is variable selection consistent. There are several other work which adopt such two-step procedures, such as the Lars-OLS hybrid (Efron et al., 2004), the relaxed Lasso (Meinshausen, 2007), the sure independence screening (Fan and Lv, 2008), the one-step sparse estimator (Zou and Li, 2008), etc. Most of the two-step procedures are computationally simple and do not require the irrepresentable condition on the design matrix, and some of them have been shown to have the oracle properties under certain conditions. However, the success of such two-step procedures depends crucially on the existence of a good initial estimator, which is not trivial to establish and also requires conditions on the design matrix especially for high-dimensional problems. For instance, Huang et al. (2006) used the univariate regression as the initial estimator in the adaptive Lasso and showed that a partial orthogonal condition is needed in order for it to satisfy the required condition in the second step.

In this paper we study several two-step procedures as well as the Ridge estimator as the initial estimator for high-dimensional problems. In Section 2 we first study under which conditions the nonnegative Garrote, adaptive Lasso and hard-thresholding procedures can turn an $\ell_\infty$-consistent estimator into a final estimator that is variable selection consistent. With some minor conditions on the penalty parameter $\lambda_n$, we show that both the nonnegative Garrote and adaptive Lasso estimators also have the oracle properties as defined in Fan and Li (2001). In Section 3 we study the conditions under which the Ridge estimator is $\ell_\infty$-consistent. The condition on the design matrix and true parameter is usually satisfied when $p_n \leq n$ and does not require the partial orthogonal condition (Huang et al., 2008) when $p_n > n$. Encouraging numerical results are provided in Section 4. Those two-step procedures with the Lasso or Ridge estimator as initial estimator are shown to have a higher success rate in terms of sparsity recovery than both the Lasso and adaptive Lasso with univariate regression as initial estimator. Results on prediction error also show that the adaptive Lasso with the Ridge initial estimator becomes more favorable when there exist stronger correlations between covariates.

## II. Two-step Procedures for Variable Selection

In the following we assume that an initial estimator $\widehat{\beta}^{init}$ can be obtained. For notational simplicity, we will use $\widehat{\beta}$ to denote the initial estimator, and also define $\Delta^* = \mathsf{diag}(\beta_1^*, \ldots, \beta_{p_n}^*)$ and $\widehat{\Delta} = \mathsf{diag}(\widehat{\beta}_1, \ldots, \widehat{\beta}_{p_n})$ respectively. We study several two-step procedures obtained using $X$, $Y$ and the initial estimator $\widehat{\beta}$.

We use $\beta_S^*$ to represent the subvector of $\beta^*$ which only contains entries $j \in S$, and it is obvious that $\beta_{S^c}^* = \mathbf{0}$. Similarly we use $X_S$ and $X_{S^c}$ to denote sub-matrices of the design matrix $X$ which only contains columns in $S$ and $S^c$, respectively. Since we are mainly interested in the situation with $p_n$ increasing, we also define $\rho_n = \min_{j \in S} |\beta_j^*|$ which is allowed to converge to zero at a relatively slow rate. Throughout the paper, we assume that $\|\beta^*\|_\infty = \max_j |\beta_j^*| < \infty$.

**Assumption 1** *Assume that the initial estimator $\widehat{\beta}$ is an $\ell_\infty$-consistent estimator of $\beta^*$, and $\|\widehat{\beta} - \beta^*\|_\infty = \max_j |\widehat{\beta}_j - \beta_j^*| = O_p(\delta_n)$ for some sequence $\delta_n \to 0$ such that $\delta_n = o(\rho_n)$.*

Although we assume that the initial estimator is a good approximation to the true parameter $\beta^*$, we do not assume that $\widehat{\beta}$ can exactly recover the sparsity pattern of $\beta^*$, since that often requires a stronger condition on the design matrix, as in the the case of the Lasso estimator. It turns out that for two-step procedures to be variable selection consistent, the $\ell_\infty$-consistent condition is sufficient. Note that similar conditions for the initial estimator have been used in earlier work (Zou, 2006; Huang et al., 2006; Yuan and Lin, 2007). It should also be obvious that in order for later procedures to separate variables in $S$ from those in $S^c$, we need to have $\rho_n$ converging to zero at a slower rate than $\delta_n$.

For any vector $\beta \in \mathbb{R}^{p_n}$, we define its support as $\mathsf{supp}(\beta) = \{j : \beta_j \neq 0\}$. A procedure is called variable selection consistent if its sequence of solutions $\widehat{\beta}_n$ as a function of sample size $n$ satisfy

$$\lim_{n \to \infty} P(\mathsf{supp}(\widehat{\beta}_n) = \mathsf{supp}(\beta^*)) = 1. \tag{2.1}$$

Furthermore, we also consider a slightly stronger property called sign consistency, which is defined by

$$\lim_{n \to \infty} P(\mathsf{sign}(\widehat{\beta}_n) = \mathsf{sign}(\beta^*)) = 1 \tag{2.2}$$

where $\mathsf{sign}(t) = -1, 0, 1$ when $t < 0$, $t = 0$ and $t > 0$ respectively. All our results about variable selection consistency trivially imply sign consistency as long as the initial estimator is $\ell_\infty$-consistent with rate faster than $\rho_n$.

### A. Nonnegative Garrote

Let $X$ and $Y$ be the design matrix and response vector, and assume that some initial estimator $\widehat{\beta}$ for the unknown parameter $\beta^*$ is given. Let $Z = X\widehat{\Delta}$, the nonnegative Garrote estimator (Breiman, 1995) $\widehat{\beta}^{NG}$ is defined as $\widehat{\beta}_j^{NG} = \widehat{\beta}_j \widehat{d}_j$ for $j = 1, \ldots, p_n$ where $\widehat{d} = (\widehat{d}_1, \ldots, \widehat{d}_{p_n})^T$ is the minimizer of

$$\frac{1}{2n}\|Y - Zd\|^2 + \lambda_n \sum_{j=1}^{p_n} d_j \tag{2.3}$$
$$d_j \geq 0 \text{ for } j = 1, \ldots, p_n. \tag{2.4}$$

Although the initial estimator for the nonnegative Garrote method was originally defined as the least squares estimator, it does not need to be so. In particular, Yuan and Lin (2007) considered a more general initial estimator for the nonnegative Garrote method with fixed $p_n$. Our result here is an extension of Yuan and Lin (2007) as we give a general sufficient condition for the nonnegative Garrote to be variable selection consistent in terms of the triple $(n, p_n, s_n)$. We start with a Lemma which is a direct consequence of the Karush-Kuhn-Tucker (KKT) condition in convex optimization.

**Lemma 2.1.** *For any $\lambda_n > 0$ and $Z = X\widehat{\Delta} = X\mathsf{diag}(\widehat{\beta}_1, \widehat{\beta}_2, \ldots, \widehat{\beta}_{p_n})$ where $\widehat{\beta}$ is some initial estimator of $\beta^*$, assume that $(Z_S^T Z_S)^{-1}$ exists. Then there exists a solution of the nonnegative Garrote that exactly recovers the sparsity pattern if and only if*

$$\left(\frac{1}{n}Z_S^T Z_S\right)^{-1} \left(\frac{1}{n}Z_S^T X_S \beta_S^* + \frac{1}{n}Z_S^T \epsilon - \lambda_n \mathbf{1}\right) > \mathbf{0} \tag{2.5}$$

$$\frac{1}{n}Z_{S^c}^T \left(I - Z_S(Z_S^T Z_S)^{-1} Z_S^T\right) \epsilon + \lambda_n Z_{S^c}^T Z_S (Z_S^T Z_S)^{-1} \mathbf{1} \leq \lambda_n \mathbf{1} \tag{2.6}$$

*where $\mathbf{0}$ and $\mathbf{1}$ are vectors composed of 0's and 1's respectively, and the inequalities hold element-wise.*

The assumption that the $s_n \times s_n$ matrix $Z_S^T Z_S$ is invertible is quite reasonable. It implies two conditions: (1) $(X_S^T X_S)^{-1}$ exists; (2) $\widehat{\beta}_j \neq 0$ for all $j \in S$. The first condition is usually needed in order to estimate $\beta_S^*$, and the second condition is satisfied as long as the initial estimator $\widehat{\beta}_S$ is element-wise close to the true parameter $\beta_S^*$ asymptotically. Furthermore, inequality (2.5) and (2.6) imply that there is no under-selection and over-selection, respectively.

We will use $\Lambda_{\min}(.)$ to denote the minimum eigenvalue operator, and in particular, we also use $\Lambda_{\min}$ to denote the lower bound of $\Lambda_{\min}(X_S^T X_S/n)$. The following result gives the conditions of the sparsity level $s_n$, the total number of predictors $p_n$ and the regularization parameter $\lambda_n$ under which the nonnegative Garotte estimator $\widehat{\beta}^{NG}$ (or $\widehat{d}$ equivalently) can correctly recover the sparsity pattern as $n \to \infty$. In other words, the nonnegative Garrote procedure is variable selection consistent when $\widehat{\beta}$ is a good initial estimator and the quantities $(n, p_n, s_n, \lambda_n, \rho_n, \delta_n)$ satisfy certain conditions.

**Theorem 2.2.** *(Nonnegative Garrote) Under Assumption 1 and further assume that*

$$\|X_{S^c}^T X_S (X_S^T X_S)^{-1}\|_\infty \leq C_{\max} < +\infty \tag{2.7}$$

$$\Lambda_{\min}\left(\frac{1}{n}X_S^T X_S\right) \geq \Lambda_{\min} > 0. \tag{2.8}$$

*Then the nonnegative Garrote estimator $\widehat{\beta}^{NG}$ is variable selection consistent, i.e.*

$$\lim_{n \to \infty} P\left(\mathsf{sign}(\widehat{\beta}^{NG}) = \mathsf{sign}(\beta^*)\right) \to 1 \tag{2.9}$$

as $n \to \infty$, if the following conditions hold:

$$\frac{\lambda_n \sqrt{s_n}}{\rho_n^2} \to 0, \quad \frac{1}{\rho_n}\sqrt{s_n \log s_n/n} \to 0, \quad \frac{\delta_n}{\lambda_n}\sqrt{\log p_n/n} \to 0. \qquad (2.10)$$

First, the irrepresentable condition for the Lasso is $\|X_{S^c}^T X_S (X_S^T X_S)^{-1} \mathsf{sign}(\beta_S^*)\|_\infty < 1$. A slightly stronger condition that does not depend on $\beta^*$ is $\|X_{S^c}^T X_S (X_S^T X_S)^{-1}\|_\infty < 1$. Here we only need to have $\|X_{S^c}^T X_S (X_S^T X_S)^{-1}\|_\infty \leq C_{\max} < \infty$ for the nonnegative Garrote if we have some good initial estimator $\widehat{\beta}$. This is mainly because

$$\|Z_{S^c}^T Z_S (Z_S^T Z_S)^{-1}\|_\infty = \left\|\widehat{\Delta}_{S^c}^T X_{S^c}^T X_S (X_S^T X_S)^{-1} \widehat{\Delta}_S^{-1}\right\|_\infty \qquad (2.11)$$

$$\leq O_p(\delta_n/\rho_n) \|X_{S^c}^T X_S (X_S^T X_S)^{-1}\|_\infty \qquad (2.12)$$

and $\delta_n = o(\rho_n)$. Also, the boundedness of $C_{\max}$ and $\Lambda_{\min}$ in equation (2.7) and (2.8) are only assumed to simplify the results and more general conditions can be obtained by allowing them converging to $\infty$ and 0 slowly. In practice, one may set the penalty parameter $\lambda_n$ proportional to $\sqrt{\log p_n/n}$. Assuming $\rho_n$ is bounded away from 0, the above conditions state that $p_n$ can increase almost as fast as $\exp(n)$, which is a well-known condition about $(p_n, n)$ for the Lasso in high-dimensional variable selection. The stringent condition on the design matrix now has been replaced by the condition that we have a good estimator $\widehat{\beta}$ such that $\max_j |\widehat{\beta}_j - \beta_j^*| = O_p(\delta_n)$.

Properties of the nonnegative Garrote estimator were studied in (Yuan and Lin, 2007) for fixed $p_n$. Although it was suspected that the nonnegative Garrote estimator might be efficient in estimation, it was only shown that $\max_j |\widehat{\beta}_j^{NG} - \beta_j^*| = O_p(\delta_n)$ for a general design matrix, that is, they only showed that $\widehat{\beta}^{NG}$ is no more better than the initial estimator $\widehat{\beta}$ in terms of estimation. In the following we show that with some additional conditions, the final nonnegative Garrote estimator is in fact efficient in estimation, i.e. it has the oracle properties (Fan and Li, 2001; Fan and Peng, 2004; Huang et al., 2006).

**Theorem 2.3.** *Let $x_i^T$ be the $i$-th row vector of $X$ (i.e. $x_i$ is the $i$-th observation), and denote $x_i^T = (x_{i(S)}^T, x_{i(S^c)}^T)$. Under assumptions in Theorem 2.2 and additionally*

$$\lambda_n \sqrt{ns_n}/\rho_n \to 0, \qquad (2.13)$$

$$n^{-1/2} \max_{1 \leq i \leq n} (x_{i(S)}^T x_{i(S)})^{1/2} \to 0, \qquad (2.14)$$

*then,*

$$\sqrt{n} w_n^{-1} v_n^T (\widehat{\beta}_S^{NG} - \beta_S^*) \to_D N(0, 1), \qquad (2.15)$$

*where $w_n^2 = \sigma^2 v_n^T \left(\frac{1}{n} X_S^T X_S\right)^{-1} v_n$ for any $s_n \times 1$ vector $v_n$ satisfying $\|v_n\|_2 \leq 1$.*

Condition 2.14 is usually satisfied if we normalize covariates and $s_n$ does not increase too fast. Condition 2.13 says $\lambda_n$ should converge to zero at a rate faster that $n^{-1/2}$ to ensure efficient estimation. In particular, if we assume $\rho_n$ is bounded away from zero, $s_n = O(1)$, $p_n = \exp(n^{1-c_1})$ and $\delta_n = n^{-1/2}$, then condition 2.10 in Theorem 2.2 together with condition 2.13 can be satisfied if we choose $\lambda_n = n^{-c_2}$ with $\frac{1}{2} < c_2 < \frac{1+c_1}{2}$.

## B. Adaptive Lasso

Given some initial estimator $\widehat{\beta}$ and define $Z = X\widehat{\Delta}$, the adaptive Lasso estimator (Zou, 2006) $\widehat{\beta}^{ALasso}$ is defined by

$$\widehat{\beta}^{ALasso} = \arg\min_{\beta} \frac{1}{2n}\|Y - X\beta\|^2 + \lambda_n \sum_{j=1}^{p_n} |\widehat{\beta}_j|^{-\gamma}|\beta_j| \qquad (2.16)$$

where $\gamma > 0$ is some tuning parameter. Considering the case $\gamma = 1$, it is easy to see that the above definition is equivalent to $\widehat{\beta}_j^{ALasso} = \widehat{\beta}_j \widehat{d}_j$ for $j = 1, \ldots, p_n$ with $\widehat{d}$ being the minimizer of

$$\widehat{d} = \arg\min_{d} \frac{1}{2n}\|Y - Zd\|^2 + \lambda_n \sum_{j=1}^{p_n} |d_j|. \qquad (2.17)$$

Zou (2006) studied properties of the adaptive Lasso for fixed $p_n$ and showed that it has the oracle properties.

The adaptive Lasso and the nonnegative Garrote, both depending on some initial estimator, are in fact closely related. It was pointed out in (Zou, 2006; Yuan and Lin, 2007) that solution of the nonnegative Garrote coincides with solution of the adaptive Lasso when additional constraints $\widehat{\beta}_j \beta_j^* \geq 0$ ($j = 1, \ldots, p_n$) are imposed. Consequently, those two methods behave very similarly when the initial estimator is of high quality. The following Lemma (Wainwright, 2006), similar to Lemma 2.1, follows from the KKT condition of the adaptive Lasso optimization problem.

**Lemma 2.4.** *For any $\lambda_n > 0$ and $Z = X\widehat{\Delta} = X\mathsf{diag}(\widehat{\beta}_1, \widehat{\beta}_2, \ldots, \widehat{\beta}_{p_n})$ where $\widehat{\beta}$ is some initial estimator of $\beta^*$, assume that $(Z_S^T Z_S)^{-1}$ exists. Then there exists a solution of adaptive Lasso that exactly recovers the sparsity pattern if and only if*

$$\left| d_S^* + \left(\frac{1}{n} Z_S^T Z_S\right)^{-1} \left(\frac{1}{n} Z_S^T \epsilon - \lambda_n \mathsf{sign}(d_S^*)\right) \right| > \mathbf{0} \qquad (2.18)$$

$$\left| Z_{S^c}^T Z_S \left(Z_S^T Z_S\right)^{-1} \left(\frac{1}{n} Z_S^T \epsilon - \lambda_n \mathsf{sign}(d_S^*)\right) - \frac{1}{n} Z_{S^c}^T \epsilon \right| \leq \lambda_n \mathbf{1} \qquad (2.19)$$

*where $\mathbf{0}$ and $\mathbf{1}$ are vectors composed of 0's and 1's, and the inequalities hold element-wise.*

The following two theorems show that under exactly the same conditions as the nonnegative Garrote, the adaptive Lasso has the oracle properties. Similar result for the adaptive Lasso has been obtained in Huang et al. (2006).

**Theorem 2.5.** *(Adaptive Lasso) Under the same conditions as in Theorem 2.2, the adaptive Lasso estimator $\widehat{\beta}^{ALasso}$ is variable selection consistent, i.e.*

$$\lim_{n\to\infty} P\left(\mathsf{sign}(\widehat{\beta}^{ALasso}) = \mathsf{sign}(\beta^*)\right) \to 1. \tag{2.20}$$

**Theorem 2.6.** *Under the same conditions as in Theorem 2.3, the adaptive Lasso estimator $\widehat{\beta}^{ALasso}$ satisfies*

$$\sqrt{n}w_n^{-1}v_n^T(\widehat{\beta}_S^{ALasso} - \beta_S^*) \to_D N(0,1), \tag{2.21}$$

*where $w_n^2 = \sigma^2 v_n^T \left(\frac{1}{n}X_S^T X_S\right)^{-1} v_n$ for any $s_n \times 1$ vector $v_n$ satisfying $\|v_n\|_2 \leq 1$.*

### C. Hard-thresholding

The hard-thresholding procedure is extremely simple and efficient. Given some initial estimator $\widehat{\beta}$ and $\lambda_n > 0$, define the hard-thresholding estimator as

$$\widehat{\beta}_j^{HT} = \begin{cases} \widehat{\beta}_j, & \text{if } |\widehat{\beta}_j| \geq \lambda_n \\ 0, & \text{if } |\widehat{\beta}_j| < \lambda_n. \end{cases} \tag{2.22}$$

Then we have the following results.

**Theorem 2.7.** *(Hard-Thresholding) Under Assumption 1 and choose $\lambda_n$ such that $\delta_n = o(\lambda_n)$ and $\lambda_n = o(\rho_n)$. Then the hard-thresholding estimator $\widehat{\beta}^{HT}$ is variable selection consistent, i.e.*

$$\lim_{n\to\infty} P\left(\mathsf{sign}(\widehat{\beta}^{HT}) = \mathsf{sign}(\beta^*)\right) \to 1. \tag{2.23}$$

Thus this simple hard-thresholding estimator can achieve variable selection consistency as well if given some good initial estimator $\widehat{\beta}$. Compared to the previous two methods, it can be directly obtained without any sophisticated optimization and has no restriction on how fast the number of variables $p_n$ and the number of relevant variables $s_n$ can grow. On the other hand, it requires that the rate of the threshold $\lambda_n$ must be greater than $\delta_n$ to ensure the variable selection consistency no matter how fast $p_n$ grows. Such an explicit relation is not needed for both the nonnegative Garrote and the Lasso, since a smaller growth rate of $p_n$ can make $\frac{\delta_n}{\lambda_n}\sqrt{\log p_n/n} \to 0$ even if $\delta_n > \lambda_n$. Hence the choice of $\lambda_n$ for the hard-thresholding procedure is more sensitive, as least from the theoretical perspective. Furthermore, it is obvious that the convergence rate of the resulting estimator $\widehat{\beta}^{HT}$ keeps the same as $\widehat{\beta}$, i.e. we have $\max_j |\widehat{\beta}_j^{HT} - \beta_j^*| = O_p(\delta_n)$. However, it is possible to apply yet another fitting method using only the subset of selected variables to obtain much better rate of convergence.

In practice, we may simply choose the hard-thresholding procedure when we know the initial estimator is $\ell_\infty$-consistent with fast convergence rate. Otherwise, the adaptive Lasso or the nonnegative Garrote might be preferred for the second step estimation and selection. We found that the latter two approaches are quite similar in terms of both theoretical properties and finite sample performance as we will see in Section 4.

## III. Initial Estimators

Clearly the success of all previous procedures crucially depends on the existence of a good initial estimator, in the sense that $\max_j |\widehat{\beta}_j - \beta_j^*| = O_p(\delta_n)$ for some sequence $\delta_n \to 0$. For $p_n$ fixed we could use the ordinary least squares (OLS) solution as the initial estimator. For $p_n$ increasing we have several choices. The simplest one is to use univariate regression (aka marginal regression), which calculates the estimator coordinate by coordinate separately, i.e. $\widehat{\beta}^{Univ} = X^T Y$. Huang et al. (2006, 2008) have used univariate regression as an initial estimator in their paper for the high-dimensional adaptive Lasso, and showed that under some partial orthogonality condition and other conditions the univariate regression estimator guarantees the zero-consistency that is closely related to the $\ell_\infty$-consistency. The partial orthogonal condition, which states that $\frac{1}{n} X_{S^c}^T X_S = O(1/\sqrt{n})$, in fact implies the irrepresentable condition asymptotically as long as $s_n$ does not grow too fast.

Another choice is to run Lasso first and use $\widehat{\beta}^{Lasso}$ as the initial estimator. Lounici (2008) studied the $\ell_\infty$ convergence rate of both the Lasso and the Dantzig selector (Candes and Tao., 2008), which requires the off-diagonal elements of $\frac{1}{n} X^T X$ to be small. Unfortunately, such a condition is quite strong and in fact implies the irrepresentable condition on the design matrix. Meinshausen and Yu (2006) showed that the Lasso estimator $\widehat{\beta}^{Lasso}$ is $\ell_2$-consistent under some sparse eigenvalue conditions. Since $\ell_2$-consistency $\|\widehat{\beta}^{Lasso} - \beta^*\|_2 = o_p(1)$ implies that $\|\widehat{\beta}^{Lasso} - \beta^*\|_\infty = O_p(\delta_n)$ for some $\delta_n \to 0$, we can use the Lasso estimator as our initial estimator. They also pointed out that the conditions under which the Lasso is $\ell_2$-consistent are not as strong as the irrepresentable condition which could fail easily even if $p_n < n$ and the design matrix is of full rank. Other works which study the $\ell_1$ or $\ell_2$-consistency of the Lasso include Bickel et al. (2007), van de Geer (2006) and Zhang and Huang (2008), which require similar sparse eigenvalue conditions on the design matrix.

We now consider another popular regression technique, the Ridge regression (Hoerl and Kennard, 1970a,b), which is more suitable for regression with correlated predictors. The Ridge estimator $\widehat{\beta}^{Ridge}$ is defined as the minimizer of the following objective (for some $\nu_n > 0$):

$$\widehat{\beta}^{Ridge} = \arg\min_\beta \frac{1}{n} \|Y - X\beta\|^2 + \nu_n \|\beta\|_2^2. \tag{3.1}$$

Our main result is that with a properly chosen regularization parameter $\nu_n$, the Ridge estimator $\widehat{\beta}^{Ridge}$ is $\ell_\infty$-consistent and thus satisfies our condition as an initial estimator. The following key assumption is needed in order to establish the $\ell_\infty$-consistent result.

**Assumption 2** *Let $\mathbf{e}_1, \ldots, \mathbf{e}_q, \mathbf{e}_{q+1}, \ldots, \mathbf{e}_{p_n}$ be the singular vectors of the symmetric matrix $\frac{1}{n} X^T X$ that corresponding to the singular values $d_1 \geq \ldots \geq d_q > d_{q+1} = \ldots = d_{p_n} = 0$ where $q$ is the rank of $\frac{1}{n} X^T X$ satisfying $q \leq \min(n, p_n)$, and let $\beta^* = \sum_{j=1}^{p_n} \theta_j \mathbf{e}_j$. Assume that $\|\sum_{j=q+1}^{p_n} \theta_j \mathbf{e}_j\|_\infty = O(\xi_n)$ with some sequence $\xi_n \to 0$.*

The requirement $\|\sum_{j=q+1}^{p_n} \theta_j \mathbf{e}_j\|_\infty = O(\xi_n)$ is obviously weaker than $\sum_{j=q+1}^{p_n} \theta_j^2 = O(\xi_n^2)$. Assumption 2 essentially says that the majority mass of $\beta^*$ belongs to the column space of $\frac{1}{n} X^T X$ asymptotically, i.e. $\beta^* \approx (\frac{1}{n} X^T X) \mathbf{b}$ for some $\mathbf{b} \in \mathbb{R}^{p_n}$ as $n \to \infty$. First, notice that the assumption is automatically satisfied when $n \leq p_n$ and $X^T X$ has full rank. However, this is not the case for the irrepresentable condition which still requires that those irrelevant predictors cannot be represented by the relevant predictors in the true model. When $p_n \gg n$ and $X^T X$ is singular, let us consider the set $\Theta = \{\theta : X\beta^* = X\theta\}$. In this case, although any $\theta \in \Theta$ is equally good in terms of predicting $Y$, there is only one true parameter $\beta^*$ among many choices. For any penalized linear method to recover the true parameter $\beta^*$, its penalty term has to favor $\beta^*$ over any other $\theta \in \Theta$. The condition in Assumption 2 can be thought as some relaxed identifiable condition for the Ridge regression to be $\ell_\infty$-consistent.

**Theorem 3.1.** *Under Assumption 2 the Ridge estimator $\widehat{\beta}^{Ridge}$ satisfies the condition $\max_j |\widehat{\beta}_j^{Ridge} - \beta_j^*| = o_p(1)$ as long as*

$$\frac{\log p_n}{n\nu_n} \to 0 \quad \text{and} \quad \frac{\nu_n \sqrt{s_n}}{d_q} \to 0. \tag{3.2}$$

*Furthermore, letting $\nu_n = (\frac{d_q^2 \log p_n}{n s_n})^{1/3}$ and if $\xi_n = O(\nu_n \sqrt{s_n}/d_q)$, we have*

$$\max_j |\widehat{\beta}_j^{Ridge} - \beta_j^*| = O_p\left(\left(\frac{\sqrt{s_n} \log p_n}{n d_q}\right)^{1/3}\right). \tag{3.3}$$

First of all, note that when $d_q$ is bounded away from 0 and $s_n = O(1)$, the result holds for $p_n = \exp(n^{1-c_1})$, $\nu_n = n^{-c_2}$ as long as $c_1 > c_2 > 0$. Such conditions can be easily satisfied for most high-dimensional linear regression problems. Notice that for the Ridge estimator to be $\ell_\infty$-consistent, there is no constraint putting on the $\rho_n$ as small coefficients do not play as important roles as in the case of variable selection. When Assumption 2 does not hold, it is easy to see that the results of Theorem 3.1 still holds for $\beta^*$'s projection $\sum_{j=1}^q \theta_j \mathbf{e}_j$. The following result shows that unlike the $\ell_\infty$-consistency, the Ridge estimator is in general not $\ell_2$-consistent with a diverging number of parameters.

**Corollary 3.2.** *The Ridge estimator $\widehat{\beta}^{Ridge}$ is in general not $\ell_2$-consistent even when $\beta^*$ is sparse and $p_n < n$.*

The main reason for the ridge estimator not being $\ell_2$-consistent is because the large number of parameters cancel out the increasing sample size. The Lasso, under certain assumptions (Meinshausen and Yu, 2006), does not suffer from such large accumulated variance due to its sparse solution. Fortunately, the two-step procedures only require the weaker $\ell_\infty$-consistency to be satisfied.

# IV. Numerical Studies

We conduct numerical experiments to evaluate finite sample properties of those two-step procedures. We consider the usage of univariate regression, OLS regression, ridge regression and the Lasso as initial estimators. These initial estimators are then processed by the nonnegative Garrote, adaptive Lasso or hard-thresholding to obtain the final estimator. In all experiments we consider the linear model $Y = X\beta^* + \epsilon$ with $\epsilon \sim \mathcal{N}(0, \sigma^2 I)$.

## A. Irrepresentable Condition and Variable Selection Consistency

First we examine how badly the irrepresentable condition will affect the success rate of those approaches. We consider an example used in (Zhao and Yu, 2007) which is to show the relationship between the probability of selecting the true sparse model and the irrepresentable condition number $\eta_\infty$ defined as:

$$\eta_\infty = 1 - \|X_{S^c}^T X_S (X_S^T X_S)^{-1} \text{sign}(\beta_S^*)\|_\infty. \tag{4.1}$$

We use the same setting as in (Zhao and Yu, 2007) by taking $n = 100$, $p = 32$ and $s = 5$, with the true sparse parameter $\beta_S^* = (7, 4, 2, 1, 1)$. The noise level $\sigma^2$ is set to 0.1 to manifest the asymptotic properties of the estimators.

We first sample a covariance matrix $\Sigma$ from $\mathsf{Wishart}(p, I_p)$, and then each sample is generated from $\mathcal{N}(0, \Sigma)$. Such a design matrix $X$ may or may not satisfy the strong irrepresentable condition (Zhao and Yu, 2007), and the degree of violation can be represented by the quantity $\eta_\infty$. When $\eta_\infty > 0$ the irrepresentable condition holds, and when $\eta_\infty < \infty$ we expect the Lasso to fail in identifying the sparsity pattern for certain cases. We generate 100 designs, and compute their corresponding $\eta_\infty$. For each design, 1000 simulations are conducted by generating the noise vector from $\mathcal{N}(0, \sigma^2 I)$. For those two-step procedures we use the Ridge regression as the initial estimator, for which the tuning parameter $\nu_n$ is automatically chosen by the generalized cross-validation (GCV). The tuning parameter $\lambda_n$ for the second step is chosen optimally over the solution path to find the correct model if possible. For Lasso we also select its optimal tuning parameter $\lambda_n^*$ by searching over the whole solution path. The advantage of using such $\lambda_n^*$ is that our variable selection results will only depend on different methods.

Figure 1 shows the percentage of correctly selected model as a function of $\eta_\infty$, and each design is shown as a dot in the plot. It is obvious that variable selection accuracy of the Lasso depends crucially on the irrepresentable condition, even for fixed $p_n$. On the other hand, results for those two-step procedures are much more accurate in terms of identifying the true model. In particular, both the nonnegative Garrote and the adaptive Lasso give almost perfect sparsity recovery for this example, with result of the hard-thresholding procedure slightly worse.

## B. High-dimensional Variable Selection Accuracy

The above example illustrates how badly the irrepresentable condition is affecting the variable selection accuracy for the Lasso. Even worse, Zhao and Yu (2007) have shown in simulation that the irrepresentable condition fails with very large probability for medium $p$ and $s$ when the design is sampled from a general Wishart distribution.

We further conduct experiment to compare the performance of different variable selection methods under a general setting. Similar to the previous example we use GCV to select $\nu_n$ for the initial estimator when applicable and use the optimal tuning parameter $\lambda_n^*$ for the second step as well as for the Lasso by search the full solution path. We let $\sigma^2 = 0.5$, $n = 50$, $p = 16, 32, 64, 128, 256, 512$ and for each $p$ we set $s = \frac{1}{16}p, \frac{2}{16}p, \ldots, \frac{15}{16}p$ unless it is greater than $n$. For each $(n, p, s)$ combination, we sample 100 times the covariance matrix $\Sigma$ from a Wishart distribution $\mathsf{Wishart}(p, I)$ and for each covariance matrix $\Sigma$ we sample every $\beta_j^*$ $(j \in S)$ uniformly from $[-2, -0.5] \cup [0.5, 2]$. For each $\Sigma$ we sample 100 times the design matrix $X$ from the multivariate normal distribution $\mathcal{N}(0, \Sigma)$. So in total there will be $100 \times 100 = 10000$ simulations for each method with the same set of $(n, p, s)$. Since we observe that results for the nonnegative Garrote and the adaptive Lasso are very similar to each other, we only report those of the adaptive Lasso. Also, we only report results for which at least one of the compared methods have success rate greater than 0.01.

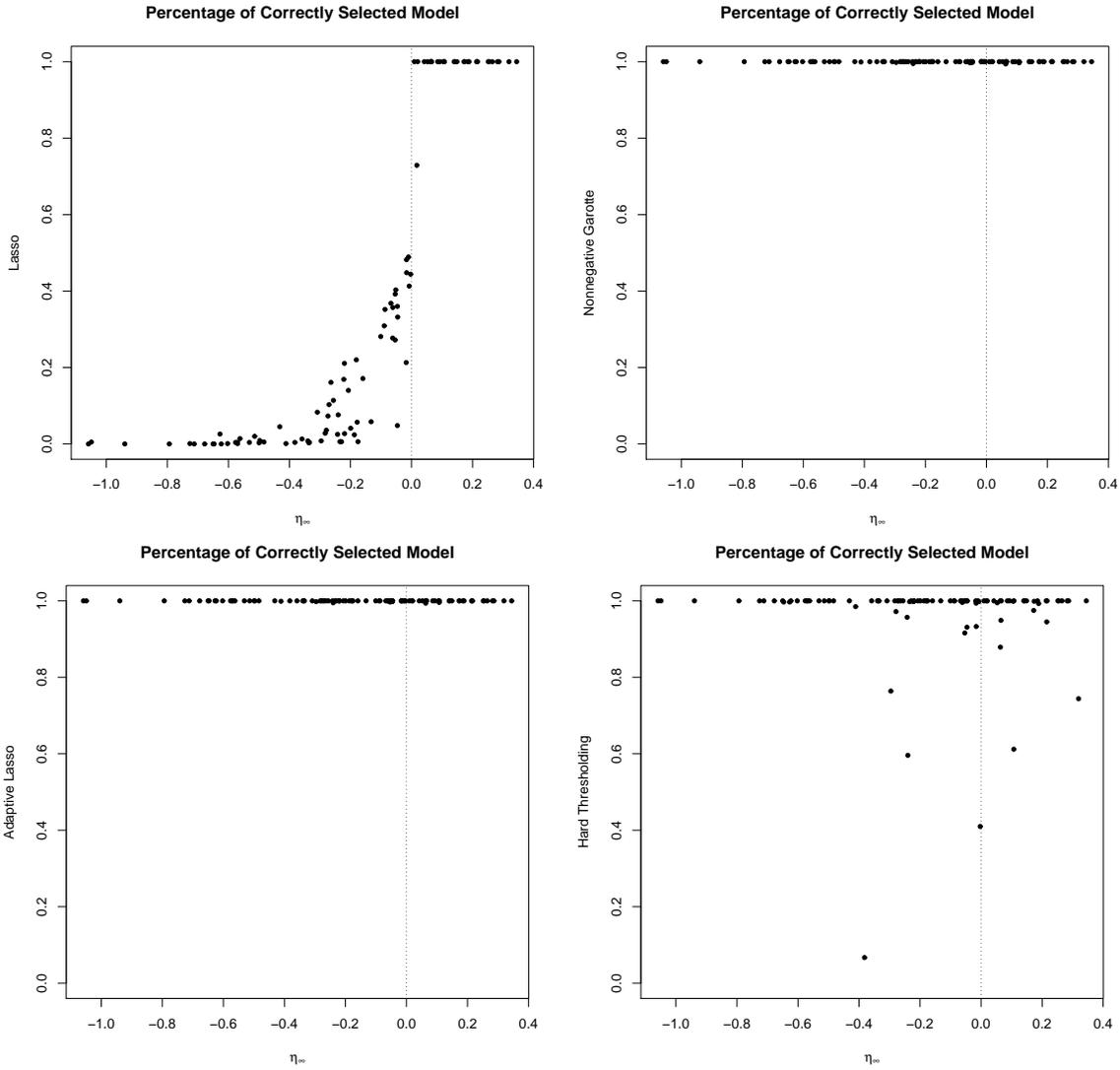

Figure 1: Example A: Percentage of correctly selected model as a function of $\eta_\infty$ for the Lasso, NG-Ridge, ALasso-Ridge and HT-Ridge: The tuning parameter $\nu_n$ for the Ridge initial estimator is chosen by GCV and the tuning parameter $\lambda_n$ is set to the optimal one by searching the full solution path.

In Table 1 the Lasso, HT-Univ and ALasso-Univ perform the worst among all the methods even for small $p$. We believe this is because of their strigent condition on the design matrix in order to achieve variable selection consistency. The two-step procedures with the Ridge initial estimator perform well especially when $s \approx p \leq n$, and those with the Lasso initial estimator performs significantly better than others when $p > n$ and $\beta^*$ is sparse.

## C. Prediction Accuracy and Variable Selection in High Dimensions

We would like to compare the following procedures: the Lasso, adaptive Lasso with univariate regression as initial estimator (ALasso-Univ), adaptive Lasso with Lasso as initial estimator (ALasso-Lasso), adaptive Lasso with Ridge as initial estimator (ALasso-Ridge), hard-thresholding with univariate regression as initial estimator (HT-Univ), hard-thresholding with Lasso as initial estimator (HT-Lasso) and hard-thresholding with Ridge as initial estimator (HT-Ridge).

To compare their prediction performance, we replicate 200 times in all the examples, and each time we generate a training dataset with 50 observations and a test dataset with 1000 observations. We use the LARS algorithm (Efron et al., 2004) to compute the Lasso and adaptive Lasso. The tuning parameter $\lambda_n$ are selected by five-fold cross validation. To measure estimation accuracy we use the relative prediction errors (RPE) defined as $E[(\widehat{y} - x^T \beta^*)^2]/\sigma^2$, and for variable selection we use the True Positive (TP) and False Positive (FP) which are defined as $TP(\beta) = \sum_{j \in S} I(\beta_j = 0)$ and $FP(\beta) = \sum_{j \notin S} I(\beta_j \neq 0)$.

*Example 1* (Auto-correlated covariance matrix). We set $p = 200$ and $\sigma = 1.5$. The covariate $x_i$ is sampled from a multivariate normal distribution with mean zero and covariance matrix $\Sigma_{j,k} = \rho^{|j-k|}$ with $\rho = 0.5, 0.75$ and $0.95$. $\beta^*$ is chosen so that there are 15 randomly located non-zero elements and the rest elements are zero. Five of the non-zero elements equal to 2.5, the second five equal to 1.5, and the last five equal to 0.5.

The auto-correlation structure of the covariance matrix in Example 1 is also used in simulations in (Tibshirani, 1996) and other Lasso related papers. It is obvious that this example is not location invariant to the variables, that is why the sparsity pattern of $\beta^*$ is randomized. Because of the high dimensionality and the degenerating feature of the covariance matrix, most of the variables are weakly correlated. Our next example has moderate to high correlations among all the variables.

Table 1: Success rate of model selection with optimally chosen $\lambda_n^*$ in the second step

| p | s | Lasso | HT- | | | | ALasso- | | | |
|---|---|---|---|---|---|---|---|---|---|---|
| | | | Univ | OLS | Ridge | Lasso | Univ | OLS | Ridge | Lasso |
| 16 | 1 | 0.9923 | 0.9787 | 0.7874 | 0.9816 | 0.9988 | 1 | 0.9257 | 0.9962 | 0.9998 |
| 16 | 3 | 0.4927 | 0.0827 | 0.7501 | 0.906 | 0.9948 | 0.4898 | 0.8596 | 0.9649 | 0.9965 |
| 16 | 5 | 0.2725 | 0.024 | 0.6721 | 0.8614 | 0.9795 | 0.2104 | 0.7919 | 0.9455 | 0.9851 |
| 16 | 7 | 0.1382 | 0 | 0.6734 | 0.7846 | 0.9529 | 0.0786 | 0.7859 | 0.8713 | 0.9557 |
| 16 | 9 | 0.0752 | 0 | 0.6392 | 0.7942 | 0.8964 | 0.0602 | 0.6875 | 0.8253 | 0.8896 |
| 16 | 11 | 0.1103 | 0.0005 | 0.657 | 0.7829 | 0.8354 | 0.051 | 0.671 | 0.7785 | 0.8139 |
| 16 | 13 | 0.173 | 0.006 | 0.7216 | 0.8336 | 0.7953 | 0.0957 | 0.7088 | 0.7896 | 0.763 |
| 16 | 15 | 0.5457 | 0.0653 | 0.7798 | 0.8443 | 0.6788 | 0.55 | 0.6863 | 0.7489 | 0.6512 |
| 32 | 2 | 0.9024 | 0.5606 | 0.6077 | 0.9714 | 0.9993 | 0.9212 | 0.8297 | 0.9952 | 0.9994 |
| 32 | 6 | 0.2027 | 0 | 0.5374 | 0.8623 | 0.9952 | 0.135 | 0.7211 | 0.9536 | 0.9974 |
| 32 | 10 | 0.0036 | 0 | 0.4137 | 0.744 | 0.9898 | 0.0021 | 0.5769 | 0.88 | 0.9924 |
| 32 | 14 | 0.0007 | 0 | 0.437 | 0.7197 | 0.9645 | 0.0005 | 0.5728 | 0.8219 | 0.9713 |
| 32 | 18 | 0.0003 | 0 | 0.4186 | 0.6984 | 0.9132 | 0 | 0.5081 | 0.7645 | 0.9095 |
| 32 | 22 | 0.0009 | 0 | 0.4536 | 0.6969 | 0.8045 | 0.0001 | 0.4829 | 0.67 | 0.7657 |
| 32 | 26 | 0.014 | 0 | 0.4827 | 0.6689 | 0.6398 | 0.0006 | 0.4346 | 0.5419 | 0.5627 |
| 32 | 30 | 0.1373 | 0.0101 | 0.5329 | 0.6951 | 0.5036 | 0.0662 | 0.3949 | 0.4669 | 0.4084 |
| 64 | 4 | 0.5752 | 0.0792 | NA | 0.8592 | 0.9993 | 0.4906 | NA | 0.9929 | 0.9999 |
| 64 | 12 | 0 | 0 | NA | 0.0662 | 0.9962 | 0 | NA | 0.3656 | 0.9994 |
| 64 | 20 | 0 | 0 | NA | 0.0086 | 0.9182 | 0 | NA | 0.0648 | 0.9257 |
| 64 | 28 | 0 | 0 | NA | 0 | 0.3995 | 0 | NA | 0.0026 | 0.4009 |
| 64 | 36 | 0 | 0 | NA | 0 | 0.0518 | 0 | NA | 0.0058 | 0.0489 |
| 64 | 44 | 0 | 0 | NA | 0.0023 | 0 | 0 | NA | 0 | 0 |
| 128 | 8 | 0.0194 | 0 | NA | 0.0048 | 1 | 0.02 | NA | 0.2633 | 1 |
| 128 | 24 | 0 | 0 | NA | 0 | 0.02 | 0 | NA | 0 | 0.02 |
| 256 | 16 | 0 | 0 | NA | 0 | 0.01 | 0 | NA | 0 | 0.01 |

*Example 2* (Constant-correlated covariance matrix). We use the same model as in Example 1 except that the covariance matrix has constant correlations, $\Sigma_{j,k} = r$ with $r = 0.3, 0.6$ and $0.85$.

The next example divides $X$ into two orthogonal blocks $X_A$ and $X_{A^c}$, so that $\Sigma_{A^c A} = 0$. Notice that when $A = S$, we have $\frac{1}{n}X_{S^c}^T X_S = O_p(1/\sqrt{n})$, which is the random version of the partial orthogonal condition for univariate estimator to be a zero-consistent initial estimator. We allow $X_A$ to be a superset of $X_S$, i.e. $A \supseteq S$.

*Example 3* (Generalized partial-orthogonal covariance matrix) We use the same model as in Example 1 except that the first 15 elements of $\beta^*$ are nonzeros and the covariance matrix has $\Sigma_{A^c A} = 0$, where $A$ includes the first $a$ columns of the $X$ and $a$ is chosen as $a = 15, 50$ and 85. All the other elements in Σ equal to constant 0.6.

Table 2: Comparing the Median RPE for Example 1 and 2 based on 200 replications†

| Example 1 | $\rho = 0.5$ | $\rho = 0.75$ | $\rho = 0.95$ |
|---|---|---|---|
| Lasso | 4.8605 (0.1783) | 4.0082 (0.1858) | 2.0365 (0.0639) |
| ALasso-Univ | 5.5771 (0.1911) | 4.6060 (0.2774) | 2.1481 (0.0658) |
| ALasso-Lasso | 4.5650 (0.3076) | 4.2502 (0.1729) | 2.5520 (0.0874) |
| ALasso-Ridge | 5.7029 (0.3080) | 4.3574 (0.2054) | 2.0348 (0.0461) |
| HT-Univ | 17.437 (0.2782) | 20.706 (0.3246) | 31.185 (0.3568) |
| HT-Lasso | 4.4008 (0.1597) | 3.5554 (0.1634) | 2.1139 (0.0795) |
| HT-Ridge | 12.122 (0.1614) | 7.6580 (0.072) | 2.4881 (0.0461) |
| Example 2 | $r = 0.3$ | $r = 0.6$ | $r = 0.85$ |
| Lasso | 3.3186 (0.1287) | 2.7097 (0.1270) | 1.8638 (0.0443) |
| ALasso-Univ | 2.9293 (0.1677) | 2.5825 (0.0872) | 1.8932 (0.0599) |
| ALasso-Lasso | 3.5494 (0.1467) | 3.1304 (0.1381) | 2.2803 (0.0655) |
| ALasso-Ridge | 3.2561 (0.2003) | 2.9072 (0.1262) | 1.8113 (1.8113) |
| HT-Univ | 66.440 (0.8814) | 68.643 (0.6295) | 35.778 (0.2719) |
| HT-Lasso | 3.0911 (0.1202) | 2.6122 (0.1538) | 1.7893 (0.0463) |
| HT-Ridge | 9.4863 (0.0707) | 5.5772 (0.0509) | 2.2337 (0.0184) |
| Example 3 | $a = 15$ | $a = 50$ | $a = 85$ |
| Lasso | 0.7355 (0.0179) | 1.3032 (0.0270) | 1.7859 (0.0399) |
| ALasso-Univ | 0.6344 (0.0166) | 1.6008 (0.0293) | 1.7467 (0.0442) |
| ALasso-Lasso | 1.1769 (0.0203) | 1.7082 (0.0422) | 2.1347 (0.0684) |
| ALasso-Ridge | 0.7217 (0.0206) | 1.3450 (0.0322) | 1.7208 (0.0361) |
| HT-Univ | 50.623 (0.4758) | 57.708 (0.4922) | 59.435 (0.5186) |
| HT-Lasso | 0.7438 (0.0162) | 1.5345 (0.0434) | 1.9722 (0.0466) |
| HT-Ridge | 3.8136 (0.0564) | 6.0480 (0.0684) | 7.3449 (0.0511) |

† The numbers in parentheses are the corresponding standard errors of RPE calculated from 200 bootstrapped sample medians.

In Table 2 Example 1, when $\rho = 0.5$ and 0.75, Lasso has better RPE than those of ALasso-Univ and ALasso-Ridge. The ALasso-Lasso and HT-lasso, which uses the Lasso as initial estimator, is also relatively good. This result is expected since the Lasso is good at dealing with situations when $s \ll p$. When $\rho = 0.95$, the ALasso-Univ and ALasso-Ridge catch up with the latter slightly better than all the other methods. In Example 2, when $r = 0.3$ and 0.6, the ALasso-Univ has better RPE than other methods. When $r = 0.85$, Alasso-ridge catches up and outperforms Lasso and other adaptive procedures. In Example 3, when $a = 15$, the partial orthogonal condition for univarate estimation is satisfied and Alasso-univ performs the best. As $a$ increases to 50, this condition is violated and Alasso-univ deteriorates faster than Lasso and other adaptive methods. In theses two cases, Lasso and

Table 3: Median number of selected variables for Example 1 and 2 based on 200 replications†

| Example 1 | $\rho = 0.5$ | | $\rho = 0.75$ | | $\rho = 0.95$ | |
|---|---|---|---|---|---|---|
| | TP | FP | TP | FP | TP | FP |
| Lasso | 11 | 18 | 11 | 24 | 9 | 26 |
| ALasso-Univ | 10 | 15 | 10 | 20 | 8.5 | 24 |
| ALasso-Lasso | 10 | 11 | 10 | 19 | 8 | 20 |
| ALasso-Ridge | 10 | 18 | 10.5 | 23.5 | 9 | 24 |
| HT-Univ | 2 | 0 | 1 | 1 | 0 | 1 |
| HT-Lasso | 11 | 17 | 11 | 23 | 8 | 17 |
| HT-Ridge | 13 | 72 | 12 | 66 | 12 | 65 |
| Example 2 | $r = 0.3$ | | $r = 0.6$ | | $r = 0.85$ | |
| | TP | FP | TP | FP | TP | FP |
| Lasso | 12 | 28 | 11 | 28 | 10 | 27 |
| ALasso-Univ | 12 | 26 | 11 | 27 | 10 | 27 |
| ALasso-Lasso | 12 | 24 | 11 | 24.5 | 9 | 22 |
| ALasso-Ridge | 12 | 25 | 11 | 25 | 10 | 24 |
| HT-Univ | 0 | 1 | 0 | 1 | 0 | 1 |
| HT-Lasso | 11 | 19 | 10 | 16 | 8 | 14 |
| HT-Ridge | 12 | 53 | 12 | 53 | 12 | 60.5 |
| Example 3 | $a = 15$ | | $a = 50$ | | $a = 85$ | |
| | TP | FP | TP | FP | TP | FP |
| Lasso | 14 | 8 | 13 | 18 | 13 | 22 |
| ALasso-Univ | 14 | 6 | 14 | 15 | 13 | 19 |
| ALasso-Lasso | 14 | 6 | 12 | 11 | 12 | 14 |
| ALasso-Ridge | 15 | 6 | 13 | 15 | 12 | 19 |
| HT-Univ | 1 | 0 | 1 | 0 | 0 | 1 |
| HT-Lasso | 14 | 2 | 13 | 10 | 12 | 17 |
| HT-Ridge | 15 | 2 | 15 | 36 | 15 | 80 |

† "TP" represents the median number of correctly selected variables, whereas "FP" represents the median number of incorrectly selected variables.

Alasso-ridge has similar RPEs. When $a$ increases to 85, Alasso-ridge outperforms all the other methods.

Hard-thresholding as another type of procedure that has different performance. HT-Univ has large RPE because of the large bias of the univariate regression. HT-Lasso however has good performance through all the cases. HT-Ridge shows up in the middle.

The variable selection results in Table 3 do not show as dramatic difference as we saw in previous examples where we choose the optimal $\lambda_n^*$ by searching the full solution path. One of the reason is that we use prediction error as the criterion to select $\lambda_n$ in the second step. Such a criterion, although could lead to good prediction accuracy, may not be ideal for the purpose of variable selection. For example, Leng et al. (2006) showed that the Lasso is not variable selection consistent in general when prediction accuracy is used as the criterion for selecting the penalty parameter. The development of an effective data-driven approach for selecting $\lambda_n$ is an interesting future research topic for variable selection.

*D. Real Data*

We study the behavior of previous methods in one real dataset to examine their predictive power. In particular, we examine the prediction accuracy of all methods as a function of sparsity level, i.e. the number of selected variables in the final model, by changing the tuning parameter $\lambda_n$. The tuning parameter $\nu_n$ for the initial estimator is chosen automatically by GCV for methods HT-Ridge, HT-Lasso, ALasso-Ridge and ALasso-Lasso.

We consider the Boston Housing data, which contains 506 records about housing values in suburbs of Boston. Each record has 13 continuous features which might be useful in describing housing price, and the response variable is the median house price. We use all 13 features as well as second order terms except for one binary feature. This results in a total of 91 predictors. In our experiments, we randomly split the data into a training set with 100 records and a test set with 406 records. We perform the random spliting 1000 times and report the average mean squared error as a function of the sparsity level of the selected model. Results are shown in Table 4. From the result we can see that the Lasso does not perform well when the sparsity level is small. This is because of the high bias for the selected variables caused by a relatively large penalty $\lambda_n$. On the other hand, those two-step procedures (except HT-Univ) do not suffer from such a problem and perform better when the sparsity level is low. As the number of selected variables increases, most methods perform reasonably well. The HT-Univ performs very poorly compared to the other two-step procedures. This is expected as the univariate estimator is not good and the hard-thresholding procedure simply cuts at a particular threshold without any data refitting.

## V. Concluding Remarks

This paper studies high-dimensional variable selection problems for linear models. In particular, we study the properties of several two-step procedures including the nonnegative Garrote, adaptive Lasso and hard-thresholding given some good initial estimator. Our results give the condition about $(n, p_n, s_n, \lambda_n)$ under which both adaptive Lasso and nonnegative Garrote can turn an $\ell_\infty$ consistent initial estimator into a final estimator that has the oracle properties as introduced by Fan and Li (2001). We then show that the Ridge estimator is $\ell_\infty$-consistent under some relaxed identifiable condition involving $\beta^*$ and $X^T X$. Such a condition is usually satisfied when $p_n \leq n$ and does not require the partial orthogonal condition needed for the univariate regression. Our simulation results show that equipped with the Lasso and Ridge estimator as initial estimators, those two-step procedures have a higher success rate in terms of sparsity recovery than the Lasso and the adaptive Lasso with the univariate regression. Results for high-dimensional estimation with correlated covariates and real data are also encouraging. Finally, it should not be difficult to extend our results to non-normal errors which have a light-tailed distribution.

Table 4: Performance of the methods as a function of sparsity level on the Boston Housing data

| sparsity | Lasso | HT- | | | ALasso- | | |
| --- | --- | --- | --- | --- | --- | --- | --- |
| | | Univ | Ridge | Lasso | Univ | Ridge | Lasso |
| 1 | 83.7634 | 67.8716 | 83.6109 | 808.557 | 60.5205 | 59.5624 | 59.29611 |
| 2 | 81.0667 | 56.7654 | 70.5460 | 347.396 | 39.0586 | 40.1758 | 45.88218 |
| 3 | 78.6111 | 73.7822 | 70.3756 | 401.654 | 34.2651 | 33.8133 | 39.78626 |
| 4 | 67.8498 | 117.6729 | 73.1718 | 261.236 | 29.8815 | 30.4250 | 35.67528 |
| 5 | 49.7597 | 181.2533 | 71.9849 | 342.206 | 27.5801 | 28.1729 | 32.65329 |
| 6 | 40.7595 | 254.4535 | 68.3086 | 271.984 | 26.2826 | 26.3106 | 29.61199 |
| 7 | 39.6022 | 337.3792 | 68.3103 | 310.752 | 25.1754 | 25.4851 | 29.16682 |
| 8 | 39.4237 | 426.7497 | 70.1900 | 387.369 | 24.2619 | 24.5633 | 29.09486 |
| 9 | 33.9314 | 521.7511 | 71.4839 | 525.376 | 23.8592 | 23.6182 | 27.54199 |
| 10 | 31.5190 | 617.0682 | 74.9066 | 401.048 | 23.5533 | 23.5084 | 27.00174 |

# VI. APPENDIX

*Proof of Lemma 2.1.*
The nonnegative Garrote is a convex optimization problem with a quadratic loss and $p_n$ linear constraints. By standard results from convex optimization we know $\widehat{d}$ is a solution of the nonnegative Garrote problem if and only if there exist $\alpha = (\alpha_1, \ldots, \alpha_{p_n})^T \geq \mathbf{0}$ such that

$$\frac{1}{n}Z^T Z\widehat{d} - \frac{1}{n}Z^T Y + \lambda_n \mathbf{1} - \alpha = \mathbf{0} \tag{6.1}$$

and $\alpha_j = 0$ if $\widehat{d}_j > 0$.

Since $\widehat{d}$ exactly recovers the sparsity pattern if and only if $\widehat{d}_{S^c} = \mathbf{0}$ and $\widehat{d}_S > \mathbf{0}$, combining these conditions with the above optimality condition we have that the nonnegative Garrote solution $\widehat{d}$ exactly recovers the sparsity pattern implies

$$\frac{1}{n}Z_S^T Z\widehat{d} - \frac{1}{n}Z_S^T Y + \lambda_n \mathbf{1} = \mathbf{0} \tag{6.2}$$

$$\frac{1}{n}Z_{S^c}^T Z\widehat{d} - \frac{1}{n}Z_{S^c}^T Y + \lambda_n \mathbf{1} \geq \mathbf{0}. \tag{6.3}$$

Since $Y = X\beta^* + \epsilon = X_S\beta_S^* + \epsilon$ and $Z\widehat{d} = Z_S\widehat{d}_S$, plugging in we have

$$\frac{1}{n}Z_S^T Z_S \widehat{d}_S - \frac{1}{n}Z_S^T X_S \beta_S^* - \frac{1}{n}Z_S^T \epsilon = -\lambda_n \mathbf{1} \tag{6.4}$$

$$\frac{1}{n}Z_{S^c}^T Z_S \widehat{d}_S - \frac{1}{n}Z_{S^c}^T X_S \beta_S^* - \frac{1}{n}Z_{S^c}^T \epsilon \geq -\lambda_n \mathbf{1}. \tag{6.5}$$

Solving the above equations we have

$$\widehat{d}_S = \left(\frac{1}{n}Z_S^T Z_S\right)^{-1}\left(\frac{1}{n}Z_S^T X_S \beta_S^* + \frac{1}{n}Z_S^T \epsilon - \lambda_n \mathbf{1}\right) \tag{6.6}$$

and

$$\frac{1}{n}Z_{S^c}^T Z_S (Z_S^T Z_S)^{-1} Z_S^T \epsilon - \lambda_n Z_{S^c}^T Z_S (Z_S^T Z_S)^{-1}\mathbf{1} - \frac{1}{n}Z_{S^c}^T \epsilon + \lambda_n \mathbf{1} \geq \mathbf{0}. \tag{6.7}$$

Now utilizing the fact that $\widehat{d}_S > \mathbf{0}$ we obtain the claimed result. □

*Proof of Theorem 2.2.*

We only need to show $\lim_{n\to\infty} P(\mathsf{supp}(\widehat{\beta}^{NG}) = \mathsf{supp}(\beta^*)) = 1$ as we have $\widehat{d} \geq \mathbf{0}$ and $\mathsf{sign}(\widehat{\beta}_S) = \mathsf{sign}(\beta_S^*)$ as $n \to \infty$ by assumption. Recall that we have diagonal matrix $\Delta^* = \mathsf{diag}(\beta_1^*, \ldots, \beta_{p_n}^*)$ and correspondingly $\widehat{\Delta} = \mathsf{diag}(\widehat{\beta}_1, \ldots, \widehat{\beta}_{p_n})$. We also use the notation $\Delta_S^*$ and $\widehat{\Delta}_S$ to denote the sub-diagonal matrices of $\Delta^*$ and $\widehat{\Delta}$ which only contains rows and columns whose indices belong to the set $S$. First, $\widehat{\Delta}_S$ is invertible with probability tending to 1 since

$$P\left(\min_{j\in S}|\widehat{\beta}_j| > 0\right) \to 1. \tag{6.8}$$

as $\delta_n = o(\rho_n)$. In the following we assume that $\widehat{\Delta}_S$ is invertible.

Define random variables $V_j = X_j^T \epsilon / n$ for $j = 1, \ldots, p_n$ and consider the events $\mathcal{A}$ and $\mathcal{B}$ given by

$$\mathcal{A} = \bigcap_{j\in S^c}\left\{|V_j| < A\sigma\sqrt{\frac{\log(p_n - s_n)}{n}}\right\} \tag{6.9}$$

$$\mathcal{B} = \bigcap_{j\in S}\left\{|V_j| < A\sigma\sqrt{\frac{\log s_n}{n}}\right\} \tag{6.10}$$

where $A$ is some constant that satisfies $A > \sqrt{2}$. By the normal error assumption we have $\sqrt{n}V_j \sim \mathcal{N}(0, \sigma^2)$, and

$$P(\mathcal{A}^c) \leq \sum_{j\in S^c} P(\sqrt{n}V_j > A\sigma\sqrt{\log(p_n - s_n)}) \tag{6.11}$$

$$\leq (p_n - s_n)P(|W| > A\sqrt{\log(p_n - s_n)}) \tag{6.12}$$

$$\leq \frac{(p_n - s_n)}{A\sqrt{\log(p_n - s_n)}} \exp\left(-\frac{A^2 \log(p_n - s_n)}{2}\right) \tag{6.13}$$

$$\leq \frac{1}{A\sqrt{\log(p_n - s_n)}} \to 0 \tag{6.14}$$

where $W$ is a standard normal variable and the last inequality is by Mill's inequality. Similarly we have

$$P(\mathcal{B}^c) \leq s_n P(|W| > A\sqrt{\log s_n}) \tag{6.15}$$

$$\leq \frac{s_n}{A\sqrt{\log s_n}} \exp\left(-\frac{A^2 \log s_n}{2}\right) \tag{6.16}$$

$$\leq \frac{1}{A\sqrt{\log s_n}} \to 0 \tag{6.17}$$

Since by our choices of events $\mathcal{A}$ and $\mathcal{B}$ we have $P(\mathcal{A} \cap \mathcal{B}) \to 1$ as $p_n > s_n \to \infty$, the following analysis will only focus on the event $\mathcal{A} \cap \mathcal{B}$. In particular, under event $\mathcal{A}$ we have the bound $\|X_{S^c}^T \epsilon/n\|_\infty < A\sigma\sqrt{\log(p_n - s_n)/n}$ and under event $\mathcal{B}$ we have the bound $\|X_S^T \epsilon/n\|_\infty < A\sigma\sqrt{\log s_n/n}$.

(1) We first show that the probability of under-selection converges to zero, and it suffices to show that

$$\widehat{d}_S = \left(\frac{1}{n}Z_S^T Z_S\right)^{-1}\left(\frac{1}{n}Z_S^T X_S \beta_S + \frac{1}{n}Z_S^T \epsilon - \lambda_n \mathbf{1}\right) \to \mathbf{1} \tag{6.18}$$

with probability 1.

Since $Z_S = X_S \widehat{\Delta}_S$, we have

$$\widehat{d}_S = \widehat{\Delta}_S^{-1} \beta_S^* + \left(\widehat{\Delta}_S^T \frac{1}{n} X_S^T X_S \widehat{\Delta}_S\right)^{-1} \widehat{\Delta}_S^T \frac{1}{n} X_S^T \epsilon - \lambda_n \left(\widehat{\Delta}_S^T \frac{1}{n} X_S^T X_S \widehat{\Delta}_S\right)^{-1} \mathbf{1}. \tag{6.19}$$

Obviously the first term converges to $\mathbf{1}$ with probability 1 at a rate $O_p(\delta_n)$ since $\delta_n = o(\rho_n)$. For the second term we have

$$\left\|\left(\widehat{\Delta}_S^T \frac{1}{n} X_S^T X_S \widehat{\Delta}_S\right)^{-1} \widehat{\Delta}_S^T \frac{1}{n} X_S^T \epsilon\right\|_\infty \leq \frac{\sqrt{s_n}}{\rho_n \Lambda_{\min}} \left\|\frac{1}{n} X_S^T \epsilon\right\|_\infty \to 0 \tag{6.20}$$

as long as $\rho_n^{-1}\sqrt{s_n \log s_n/n} \to 0$.

For the last term we have

$$\left\|\lambda_n \left(\widehat{\Delta}_S^T \frac{1}{n} X_S^T X_S \widehat{\Delta}_S\right)^{-1} \mathbf{1}\right\|_\infty \leq \frac{\lambda_n \sqrt{s_n}}{\Lambda_{\min}} \left\|\widehat{\Delta}_S^{-1}\right\|_\infty^2 = O_p\left(\frac{\lambda_n \sqrt{s_n}}{\rho_n^2}\right). \tag{6.21}$$

Combining three terms together we have $\widehat{d}_S \to \mathbf{1}$ with probability 1 if $\lambda_n = o(\rho_n^2/\sqrt{s_n})$ and $\rho_n^{-1}\sqrt{s_n \log s_n/n} \to 0$.

(2) We show that the probability of over-selection converges to 0 as well. First, Define

$$W = \frac{1}{n} Z_{S^c}^T \left(I - Z_S(Z_S^T Z_S)^{-1} Z_S^T\right) \epsilon + \lambda_n Z_{S^c}^T Z_S (Z_S^T Z_S)^{-1} \mathbf{1} \tag{6.22}$$

and there is no over-selection if $\max_{j \in S^c} W_j \leq \lambda_n$, which is further implied by the event $\|W\|_\infty \leq \lambda_n$. We have

$$\|W\|_\infty = \left\| \frac{1}{n} Z_{S^c}^T \left( I - Z_S (Z_S^T Z_S)^{-1} Z_S^T \right) \epsilon + \lambda_n Z_{S^c}^T Z_S (Z_S^T Z_S)^{-1} \mathbf{1} \right\|_\infty \tag{6.23}$$

$$\leq \left\| \frac{1}{n} Z_{S^c}^T \epsilon \right\|_\infty + \left\| \frac{1}{n} Z_{S^c}^T Z_S (Z_S^T Z_S)^{-1} Z_S^T \epsilon \right\|_\infty + \lambda_n \| Z_{S^c}^T Z_S (Z_S^T Z_S)^{-1} \mathbf{1} \|_\infty \tag{6.24}$$

$$= O_p(\delta_n) \left\| \frac{1}{n} X_{S^c}^T \epsilon \right\|_\infty + O_p(\delta_n) \left\| \frac{1}{n} X_S^T \epsilon \right\|_\infty + O_p\left( \frac{\lambda_n \delta_n}{\rho_n} \right) \tag{6.25}$$

Thus on the event $\mathcal{A} \cap \mathcal{B}$, we have $\|W\|_\infty \leq \lambda_n$ as $n \to \infty$ as long as $\frac{\delta_n}{\lambda_n} \sqrt{\log p_n / n} \to 0$. The result now follows by combining (1) and (2). □

*Proof of Theorem 2.3.*

This theorem can be verified in a similar way as in the proof of Theorem 2 of (Huang et al., 2008). By Theorem 2.2, $P(\widehat{d}_{S^c} = \mathbf{0}) \to 1$ and $P(\widehat{d}_S \neq \mathbf{0}) \to 1$, then, the KKT condition implies

$$\left( \frac{1}{n} Z_S^T Z_S \right) \widehat{d}_S - \frac{1}{n} Z_S^T Y = -\lambda_n \mathbf{1}. \tag{6.26}$$

Plug in $Y = X_S \beta_S^* + \epsilon$, $Z_S = X_S \widehat{\Delta}_S$ and $\widehat{d}_S = (\widehat{\Delta}_S)^{-1} \widehat{\beta}_S^{NG}$, we have

$$\left( \frac{1}{n} \widehat{\Delta}_S X_S^T X_S \right) \left( \widehat{\beta}_S^{NG} - \beta_S^* \right) = \frac{1}{n} \widehat{\Delta}_S X_S^T \epsilon - \lambda_n \mathbf{1}, \tag{6.27}$$

then

$$\sqrt{n} v_n^T \left( \widehat{\beta}_S^{NG} - \beta_S^* \right) = n^{-1/2} v_n^T \left( \frac{1}{n} X_S^T X_S \right)^{-1} X_S^T \epsilon - \sqrt{n} \lambda_n v_n^T \left( \frac{1}{n} X_S^T X_S \right)^{-1} (\widehat{\Delta}_S)^{-1} \mathbf{1}. \tag{6.28}$$

Since

$$\left| \sqrt{n} \lambda_n v_n^T \left( \frac{1}{n} X_S^T X_S \right)^{-1} (\widehat{\Delta}_S)^{-1} \mathbf{1} \right| \leq \sqrt{n} \lambda_n \| v_n \|_2 \left\| \left( \frac{1}{n} X_S^T X_S \right)^{-1} (\widehat{\Delta}_S)^{-1} \mathbf{1} \right\|_2 \tag{6.29}$$

$$\leq \sqrt{n s_n} \lambda_n \Lambda_{\min}^{-1} \left( \Lambda_{\min}(\widehat{\Delta}_S) \right)^{-1} \tag{6.30}$$

$$\leq \sqrt{n s_n} \lambda_n \Lambda_{\min}^{-1} \rho_n^{-1} (1 + o_p(1)), \tag{6.31}$$

then, under condition (2.13), we have

$$\sqrt{n} w_n^{-1} v_n^T \left( \widehat{\beta}_S^{NG} - \beta_S^* \right) = n^{-1/2} w_n^{-1} v_n^T \left( \frac{1}{n} X_S^T X_S \right)^{-1} X_S^T \epsilon + o_p(1). \tag{6.32}$$

Next, we verify the conditions for Linderberg-Feller central limit theorem. Let

$$V_i = n^{-1/2} w_n^{-1} v_n^T \left( \frac{1}{n} X_S^T X_S \right)^{-1} x_{i(S)}, \tag{6.33}$$

and $W_i = V_i \epsilon_i$, then it is easy to show that

$$\text{Var}\left(\sum_{i=1}^n W_i\right) = \sigma^2 \sum_{i=1}^n V_i^2 = 1. \tag{6.34}$$

On the other hand,

$$\sum_{i=1}^n \mathbb{E}\left[W_i^2 1(|W_i| > \delta)\right] = \sigma^2 \sum_{i=1}^n V_i^2 \mathbb{E}\left[\epsilon_i^2 1(|V_i \epsilon_i| > \delta)\right] \leq \max_{1 \leq i \leq n} \mathbb{E}\left[\epsilon_i^2 1(|V_i \epsilon_i| > \delta)\right], \tag{6.35}$$

then it is enough to show that

$$\max_{1 \leq i \leq n} \mathbb{E}\left[\epsilon_i^2 1(|V_i \epsilon_i| > \delta)\right] \to 0, \tag{6.36}$$

or equivalently,

$$\max_{1 \leq i \leq n} |V_i| = n^{-1/2} w_n^{-1} \max_{1 \leq i \leq n} \left| v_n^T \left(\frac{1}{n} X_S^T X_S\right)^{-1} x_{i(S)} \right| \to 0. \tag{6.37}$$

Since

$$\left| v_n^T \left(\frac{1}{n} X_S^T X_S\right)^{-1} x_{i(S)} \right| \leq \left( v_n^T \left(\frac{1}{n} X_S^T X_S\right)^{-1} v_n \right)^{1/2} \left( x_{i(S)}^T \left(\frac{1}{n} X_S^T X_S\right)^{-1} x_{i(S)} \right)^{1/2} \tag{6.38}$$

$$\leq \sigma^{-1} w_n \Lambda_{min}^{-1/2} \left( x_{i(S)}^T x_{i(S)} \right)^{1/2}, \tag{6.39}$$

then under assumption (2.14), (6.37) follows. This finishes the proof. □

*Proof of Lemma 2.4.*

By assumption we have $\widehat{\beta}_S \neq \mathbf{0}$ and thus $\widehat{\beta}^{ALasso}$ exactly recovers the sparsity pattern if and only if $\widehat{d}$ does so. By the KKT condition, $\widehat{d}$ is a solution if and only if there exists a subgradient $\widehat{z} \in \partial \ell_1(\widehat{d})$ such that

$$\frac{1}{n} Z^T Z \widehat{d} - \frac{1}{n} Z^T Y + \lambda_n \widehat{z} = \mathbf{0} \tag{6.40}$$

where $\widehat{z}_j = \text{sign}(\widehat{d}_j)$ for $\widehat{d}_j \neq 0$ and $|\widehat{z}_j| \leq 1$ otherwise. Then it follows that $\widehat{d}$ (and thus $\widehat{\beta}^{ALasso}$) exactly recovers the sparsity pattern if and only if $\widehat{d}_{S^c} = \mathbf{0}$, $\widehat{d}_S \neq \mathbf{0}$, $|\widehat{z}_{S^c}| \leq \mathbf{1}$ and $\widehat{z}_S = \text{sign}(d_S^*)$.

Combining these conditions with the above optimality condition we have that the adaptive Lasso solution $\widehat{\beta}^{ALasso}$ recovers the sparsity pattern implies

$$\frac{1}{n} Z_S^T Z \widehat{d} - \frac{1}{n} Z_S^T Y + \lambda_n \widehat{z}_S = \mathbf{0} \tag{6.41}$$

$$\frac{1}{n} Z_{S^c}^T Z \widehat{d} - \frac{1}{n} Z_{S^c}^T Y + \lambda_n \widehat{z}_{S^c} = \mathbf{0}. \tag{6.42}$$

Since $Y = Zd^* + \epsilon = Z_S d_S^* + \epsilon$ and $Z\widehat{d} = Z_S \widehat{d}_S$, plugging in we have

$$\frac{1}{n} Z_S^T Z_S \widehat{d}_S - \frac{1}{n} Z_S^T Z_S d_S^* - \frac{1}{n} Z_S^T \epsilon = -\lambda_n \text{sign}(d_S^*) \quad (6.43)$$

$$\frac{1}{n} Z_{S^c}^T Z_S \widehat{d}_S - \frac{1}{n} Z_{S^c}^T Z_S d_S^* - \frac{1}{n} Z_{S^c}^T \epsilon = -\lambda_n \widehat{z}_{S^c}. \quad (6.44)$$

Solving the above equations we have

$$\widehat{d}_S = d_S^* + \left(\frac{1}{n} Z_S^T Z_S\right)^{-1} \left(\frac{1}{n} Z_S^T \epsilon - \lambda_n \text{sign}(d_S^*)\right) \quad (6.45)$$

$$-\lambda_n \widehat{z}_{S^c} = Z_{S^c}^T Z_S \left(Z_S^T Z_S\right)^{-1} \left(\frac{1}{n} Z_S^T \epsilon - \lambda_n \text{sign}(d_S^*)\right) - \frac{1}{n} Z_{S^c}^T \epsilon \quad (6.46)$$

and the result follows since $|\widehat{d}_S| > \mathbf{0}$ and $|\widehat{z}_{S^c}| \leq \mathbf{1}$. □

*Proof of Theorem 2.5.*

The proof is similar to that of Theorem 2.2. Without loss of generality, assume that $\widehat{\Delta}$ is invertible and define events $\mathcal{A}$ and $\mathcal{B}$ as before. We only need to consider the situation when $\mathcal{A} \cap \mathcal{B}$ is true.

(1) We have $d_S^* \to \mathbf{1}$ since $\delta_n = o(\rho_n)$. As in Theorem 2.2, we have

$$\left\| \left(\frac{1}{n} Z_S^T Z_S\right)^{-1} \frac{1}{n} Z_S^T \epsilon \right\|_\infty \leq \frac{\sqrt{s_n}}{\rho_n \Lambda_{\min}} \left\| \frac{1}{n} X_S^T \epsilon \right\|_\infty \to 0 \quad (6.47)$$

as long as $\rho_n^{-1} \sqrt{s_n \log s_n / n} \to 0$. Also we have

$$\left\| \lambda_n \left(\frac{1}{n} Z_S^T Z_S\right)^{-1} \text{sign}(d_S^*) \right\|_\infty \leq \frac{\lambda_n \sqrt{s_n}}{\Lambda_{\min}} \left\| \widehat{\Delta}_S^{-1} \right\|_\infty^2 = O_p\left(\frac{\lambda_n \sqrt{s_n}}{\rho_n^2}\right) \quad (6.48)$$

Thus we have if $\lambda_n = o(\rho_n^2 / \sqrt{s_n})$ and $\rho_n^{-1} \sqrt{s_n \log s_n / n} \to 0$.

(2) Define $W = \frac{1}{n} Z_{S^c}^T \left(I - Z_S(Z_S^T Z_S)^{-1} Z_S^T\right) \epsilon + \lambda_n Z_{S^c}^T Z_S (Z_S^T Z_S)^{-1} \text{sign}(d_S^*)$ which is the same as the random vector $W$ in the proof of Theorem 2.2 except that $\mathbf{1}$ is replaced by $\text{sign}(d_S^*)$. Thus we have $\|W\|_\infty \leq \lambda_n$ if $\delta_n \sqrt{\log p_n / n} \to 0$. □

*Proof of Theorem 2.6.*

By Theorem 2.5, we have $P(\widehat{d}_{S^c} = \mathbf{0}) \to 1$ and $P(\widehat{d}_S \neq \mathbf{0}) \to 1$. Then the KKT condition implies

$$\left(\frac{1}{n} \widehat{\Delta}_S X_S^T X_S\right) \left(\widehat{\beta}_S^{ALasso} - \beta_S^*\right) = \frac{1}{n} \widehat{\Delta}_S X_S^T \epsilon - \lambda_n \text{sign}(\widehat{d}_S), \quad (6.49)$$

and the rest follows exactly as the proof of Theorem 2.3. □

*Proof of Theorem 2.7.*

For all $j \notin S$, i.e. $j$ such that $\beta_j^* = 0$, we have

$$P\left(\max_{j \notin S} |\widehat{\beta}_j| \geq \lambda_n\right) = P\left(\max_{j \notin S} |\widehat{\beta}_j/\delta_n| \geq \lambda_n/\delta_n\right) \to 0$$

since $\delta_n = o(\lambda_n)$. By the hard-thresholding rule, we have $P(\widehat{\beta}_{S^c}^{HT} = \mathbf{0}) \to 1$.

For all $j \in S$, i.e. $j$ such that $\beta_j^* \neq 0$, we have

$$P\left(\inf_{j \in S} |\widehat{\beta}_j| > \lambda_n\right) \geq P\left(\inf_{j \in S}(|\beta_j^*| - |\widehat{\beta}_j - \beta_j^*|) > \lambda_n\right) \geq P\left(\rho_n - \max_{j \in S} |\widehat{\beta}_j - \beta_j^*| > \lambda_n\right)$$

since $\inf_{j \in S} |\widehat{\beta}_j| \geq \rho_n - \max_{j \in S} |\widehat{\beta}_j - \beta_j^*|$. The right hand side converges to 1 as long as $\lambda_n = o(\rho_n)$. As a result, we have $P(\widehat{\beta}_S^{HT} = \widehat{\beta}_S) = 1$. □

*Proof of Theorem 3.1.*

First, notice that $(\widehat{\beta}^{Ridge} - \beta^*)$ is a random vector which follows a multivariate normal distribution with mean

$$-\nu_n \left(\frac{1}{n} X^T X + \nu_n I\right)^{-1} \beta^* \tag{6.50}$$

and covariance matrix

$$\text{Var}\left(\widehat{\beta}^{Ridge} - \beta^*\right) = \frac{\sigma^2}{n^2} \left(\frac{1}{n} X^T X + \nu_n I\right)^{-1} X^T X \left(\frac{1}{n} X^T X + \nu_n I\right)^{-1} \tag{6.51}$$

$$= \frac{\sigma^2}{n} \left(\left(\frac{1}{n} X^T X + \nu_n I\right)^{-1} - \nu_n \left(\frac{1}{n} X^T X + \nu_n I\right)^{-2}\right). \tag{6.52}$$

Let $\mathbf{m}$ be the mean vector and $\mathbf{C}$ be the covariance matrix of $(\widehat{\beta}^{Ridge} - \beta^*)$ respectively, and define $\bar{m} = \max_j |m_j|$ and $\bar{C} = \max_j C_{jj}$ to be the uniform upper bound of the individual bias and variance.

Define event $\mathcal{E}$ to be

$$\mathcal{E} = \bigcap_{j=1}^{p_n} \left\{|\widehat{\beta}_j^{Ridge} - \beta_j^*| \leq \sqrt{2\bar{C} \log p_n} + \bar{m}\right\}, \tag{6.53}$$

then we have

$$P(\mathcal{E}^c) \leq \sum_{j=1}^{p_n} P\left(|\widehat{\beta}_j^{Ridge} - \beta_j^*| > \sqrt{2\bar{C} \log p_n} + \bar{m}\right) \tag{6.54}$$

$$\leq \sum_{j=1}^{p_n} P\left(|\widehat{\beta}_j^{Ridge} - \beta_j^* - m_j| > \sqrt{2 C_{jj} \log p_n}\right) \tag{6.55}$$

$$= p_n P\left(|Z| > \sqrt{2 \log p_n}\right) \tag{6.56}$$

$$\leq \frac{p_n}{\sqrt{2 \log p_n}} \exp\left(-\log p_n\right) \to 0. \tag{6.57}$$

where $Z \sim \mathcal{N}(0,1)$ is a standard normal random variable. So we only need to consider the situation on the event $\mathcal{E}$. In other words, we need to bound the quantity $\sqrt{2\bar{C} \log p_n} + \bar{m}$.

We first compute $\bar{C}$. Define $D_{\max}(\mathbf{C})$ to be the operator which returns the maximum diagonal element of $\mathbf{C}$, and recall that $\Lambda_{\max}(\mathbf{C})$ is the maximum eigenvalue of matrix $\mathbf{C}$, we have

$$\bar{C} = \frac{\sigma^2}{n} D_{\max}\left(\left(\frac{1}{n}X^T X + \nu_n I\right)^{-1} - \nu_n \left(\frac{1}{n}X^T X + \nu_n I\right)^{-2}\right) \tag{6.58}$$

$$\leq \frac{\sigma^2}{n} D_{\max}\left(\left(\frac{1}{n}X^T X + \nu_n I\right)^{-1}\right) \tag{6.59}$$

$$\leq \frac{\sigma^2}{n} \Lambda_{\max}\left(\left(\frac{1}{n}X^T X + \nu_n I\right)^{-1}\right) \tag{6.60}$$

$$\leq \frac{\sigma^2}{n\nu_n}. \tag{6.61}$$

Next we bound $\bar{m}$. Since we have

$$\frac{1}{n}X^T X = UDU^T \tag{6.62}$$

where $U \in \mathbb{R}^{p \times p}$ is an orthogonal matrix and $D = \mathsf{diag}(d_1, d_2, \ldots, d_q, 0, \ldots, 0)$ is a diagonal matrix with $d_1 \geq d_2 \geq \ldots \geq d_q > 0$. Note that $q \leq n$ since $D$ has at most $n$ nonzero elements, and we also use $D^-$ to represent the pseudo-inverse of $D$. Let columns of $U$ be $\mathbf{e}_1, \ldots, \mathbf{e}_{p_n}$. By Assumption 2 we have $\beta^* = (\frac{1}{n}X^T X)\mathbf{b} + \sum_{j=q+1}^{p_n} \theta_j \mathbf{e}_j$ for some vector $\mathbf{b} \in \mathbb{R}^p$ and $\|\sum_{j=q+1}^{p_n} \theta_j \mathbf{e}_j\|_\infty = O(\xi_n)$.

Thus we have

$$\bar{m} = \left\|-\nu_n \left(\frac{1}{n}X^T X + \nu_n I\right)^{-1} \beta^*\right\|_\infty \tag{6.63}$$

$$= \left\|\nu_n U (D + \nu_n I)^{-1} U^T \beta^*\right\|_\infty \tag{6.64}$$

$$\leq \left\|\nu_n U (D + \nu_n I)^{-1} U^T \frac{1}{n}X^T X \mathbf{b}\right\|_\infty + \left\|\nu_n U (D + \nu_n I)^{-1} U^T \left(\sum_{j=q+1}^{p_n} \theta_j \mathbf{e}_j\right)\right\|_\infty \tag{6.65}$$

$$= \left\|\nu_n U (D + \nu_n I)^{-1} DD^- DU^T \mathbf{b}\right\|_\infty + \left\|\sum_{j=q+1}^{p_n} \theta_j \mathbf{e}_j\right\|_\infty \tag{6.66}$$

$$\leq \left\|\nu_n U (D + \nu_n I)^{-1} DD^- DU^T \mathbf{b}\right\|_2 + O(\xi_n) \tag{6.67}$$

$$\leq \Lambda_{\max}\left(\nu_n (D + \nu_n I)^{-1} D\right) \left\|D^- DU^T \mathbf{b}\right\|_2 + O(\xi_n) \tag{6.68}$$

$$= \frac{\nu_n d_1}{\nu_n + d_1} \left\|D^- DU^T \mathbf{b}\right\|_2 + O(\xi_n) \tag{6.69}$$

$$\leq O\left(\frac{\nu_n \sqrt{s_n}}{d_q} + \xi_n\right). \tag{6.70}$$

The last inequality comes from the fact that

$$\|D^- DU^T \mathbf{b}\|_2 = \|D^- U^T \beta^* - D^-[0,\ldots,0,\theta_{q+1},\ldots,\theta_{p_n}]^T\|_2 \leq \frac{1}{d_q} O(\sqrt{s_n})$$

since the true parameter $\beta^*$ is assumed to be sparse with only $s_n$ number of nonzero elements. Combining the above steps we get that on the event $\mathcal{E}$, we have

$$\left\|\widehat{\beta}^{Ridge} - \beta^*\right\|_\infty \leq \sqrt{2\bar{C}\log p_n} + \bar{m} \tag{6.71}$$

$$\leq \sqrt{\frac{4\sigma^2}{n\nu_n}\log p_n} + O\left(\frac{\nu_n\sqrt{s_n}}{d_q} + \xi_n\right) \tag{6.72}$$

$$\to 0 \tag{6.73}$$

as long as $\frac{\nu_n\sqrt{s_n}}{d_q} \to 0$ and $\frac{\log p_n}{n\nu_n} \to 0$. Furthermore, if $\xi_n \to 0$ sufficiently fast, we have

$$\|\widehat{\beta}^{Ridge} - \beta^*\|_\infty = O_p\left(\left(\frac{\sqrt{s_n}\log p_n}{nd_q}\right)^{1/3}\right) \tag{6.74}$$

by setting $\nu_n = (\frac{d_q^2 \log p_n}{n s_n})^{1/3}$. □

*Proof of Corollary 3.2.*

We only need to consider the special case of orthogonal design where $\frac{1}{n}X^T X = I_{p_n}$. In order to have the orthogonal design, we need to have $p_n <= n$. Suppose $p_n = n/2$ and the design is orthogonal such that $\frac{1}{n}X^T X = I_{p_n}$.

Then we have

$$\widehat{\beta}^{Ridge} - \beta^* = -\frac{\nu_n}{1+\nu_n}\beta^* + \frac{1}{1+\nu_n}\frac{1}{n}X^T\epsilon. \tag{6.75}$$

As a result, in order for $\widehat{\beta}^{Ridge}$ to be an $\ell_2$-consistent estimator of $\beta^*$, both the first and second term need to disappear. The first term goes to 0 for arbitrary $\beta^*$ only if $\lambda_n \to 0$, and in this case we need $\|X^T\epsilon/n\|_2 = o_p(1)$ to ensure the $\ell_2$-consistency. However, we have $\mathbb{E}(X^T\epsilon/n) = \mathbf{0}$ and $\text{Var}(X^T\epsilon/n) = \frac{\sigma^2}{n}I_{p_n}$. Consequently, $\|X^T\epsilon/n\|_2 = O_p(1)$ since $p_n = n/2$ and the proof is completed. □

# References

BICKEL, P. J., RITOV, Y. and TSYBAKOV, A. (2007). Simultaneous analysis of lasso and dantzig selector. *Technical report, U.C.Berkeley* .

BREIMAN, L. (1995). Better subset regression using the nonnegative garrote. *Technometrics* **37** 373–384.

CANDES, E. J. and TAO., T. (2008). The dantzig selector: statistical estimation when p is much larger than n. *The Annals of Statistics* **to appear**.

EFRON, B., HASTIE, T., JOHNSTONE, I. and TIBSHIRANI, R. (2004). Least angle regression. *The Annals of Statistics* **32** 407–499.

FAN, J. and LI, R. (2001). Variable selection via nonconcave penalized likelihood and its oracle properties. *Journal of the American Statistical Association* **96** 1348–1360.

FAN, J. and LV, J. (2008). Sure independence screening for ultra-high dimensional feature space. *Journal of the Royal Statistical Society, Series B, Methodological* **70** 849–911.

FAN, J. and PENG, H. (2004). Nonconcave penalized likelihood with a diverging number of parameters. *The Annals of Statistics* **32** 928–961.

FRANK, I. E. and FRIEDMAN, J. H. (1993). A statistical view of some chemometrics regression tools (with discussion). *Technometrics* **35** 109–148.

FU, W. and KNIGHT, K. (2000). Asymptotics for lasso type estimators. *The Annals of Statistics* **28** 1356–1378.

GREENSHTEIN, E. and RITOV, Y. (2004). Persistency in high dimensional linear predictor-selection and the virtue of over-parametrization. *Journal of Bernoulli* **10** 971–988.

HOERL, A. and KENNARD, R. (1970a). Ridge regression: Applications to nonorthogonal problems. *Technometrics* **12** 69–82.

HOERL, A. and KENNARD, R. (1970b). Ridge regression: Biased estimation for nonorthogonal problems. *Technometrics* **12** 55–67.

HUANG, J., HOROWITZ, J. L. and MA, S. (2008). Asymptotic properties of bridge estimators in sparse high-dimensional regression models. *Annals of Statistics* **36** 587–613.

HUANG, J., MA, S. and ZHANG, C. H. (2006). Adaptive lasso for sparse high-dimensional regression models. *Technical Report No. 374, Department of Statistics and Actuarial Science, University of Iowa* .

LENG, C., LIN, Y. and WAHBA, G. (2006). A note on the lasso for sparse high-dimensional regression models. *Statistical Sinica* To appear.

LOUNICI, K. (2008). Sup-norm convergence rate and sign concentration property of lasso and dantzig estimators. *Electronic Journal of Statistics* **2** 90–102.

MEINSHAUSEN, N. (2007). Relaxed lasso. *Computational Statistics and Data Analysis* **52** 374–393.

MEINSHAUSEN, N. and BÜHLMANN, P. (2006). High dimensional graphs and variable selection with the lasso. *The Annals of Statistics* **34** 1436–1462.

MEINSHAUSEN, N. and YU, B. (2006). Lasso-type recovery of sparse representations for high-dimensional data. Tech. Rep. 720, Department of Statistics, UC Berkeley.

OSBORNE, M. R., PRESNELL, B. and TURLACH, B. A. (2000). On the lasso and its dual. *Journal of Computational and Graphical Statistics* **9** 319–337.

Tibshirani, R. (1996). Regression shrinkage and selection via the lasso. *Journal of the Royal Statistical Society, Series B, Methodological* **58** 267–288.

van de Geer, S. (2006). High-dimensional generalized linear models and the lasso. *To appear in the Annals of Statistics*.

Wainwright, M. J. (2006). Sharp thresholds for high-dimensional and noisy sparsity recovery using $\ell_1$-constrained quadratic programs. In *Proc. Allerton Conference on Communication, Control and Computing*.

Yuan, M. and Lin, Y. (2007). On the non-negative garrote estimator. *Journal of the Royal Statistical Society, Series B, Methodological* **69** 143–161.

Zhang, C.-H. and Huang, J. (2008). The sparsity and bias of the lasso selection in high-dimensional linear regression. *The Annals of Statistics* **36** 1567–1594.

Zhao, P. and Yu, B. (2007). On model selection consistency of lasso. *J. of Mach. Learn. Res.* **7** 2541–2567.

Zou, H. (2006). The adaptive lasso and its oracle properties. *Journal of the American Statistical Association* **101** 1418–1429.

Zou, H. and Li, R. (2008). One-step sparse estimates in nonconcave penalized likelihood models. *The Annals of Statistics* **36** 1509–1533.